\documentclass[11pt]{article}

\usepackage{amsmath,amssymb,amsfonts,graphicx}
\usepackage[latin1]{inputenc}

\newtheorem{thm}{Theorem}[section]
\newtheorem{prop}[thm]{Proposition}
\newtheorem{lem}[thm]{Lemma}
\newtheorem{df}[thm]{Definition}

\newtheorem{rem}[thm]{Remark}

\newtheorem{cor}[thm]{Corollary}

\def\be#1 {\begin{equation} \label{#1}}
\newcommand{\ee}{\end{equation}}
\def\dem {\noindent {\bf Proof : }}

\newcommand{\mb}{\medskip\noindent}
\newcommand{\gb}{\bigskip\noindent}
\newcommand{\R}{\mathbb R}

\newcommand{\F}{\mathcal F}
\newcommand{\A}{\mathcal A}
\newcommand{\M}{\mathcal M}

\newcommand{\diam}{diam}
\newcommand{\aver}[1]{-\hskip-0.46cm\int_{#1}}

\newcommand{\ind}{ {\bf 1} }

\DeclareMathOperator{\Lip}{Lip}

\def\sqw{\hbox{\rlap{\leavevmode\raise.3ex\hbox{$\sqcap$}}$%
\sqcup$}}
\def\findem{\ifmmode\sqw\else{\ifhmode\unskip\fi\nobreak\hfil
\penalty50\hskip1em\null\nobreak\hfil\sqw
\parfillskip=0pt\finalhyphendemerits=0\endgraf}\fi}

\title{New Calder\'on-Zygmund decompositions}
\date{July 15, 2009}

\begin {document}
\author{ N. Badr\\Institut Camille Jordan\\Universit\'e Claude Bernard Lyon 1 \\ 43 boulevard du 11 Novembre 1918\\ F-69622 Villeurbanne Cedex\\ badr@math.univ-lyon1.fr \and
F. Bernicot\\ Laboratoire de Math\'ematiques \\
Universit\'e Paris-Sud XI \\F-91405 Orsay
Cedex\\frederic.bernicot@math.u-psud.fr \\}
\maketitle

\begin{abstract} 
\mb

We state a new Calder\'on-Zygmund decomposition for Sobolev spaces on a doubling Riemannian manifold. Our hypotheses are weaker than those of the already known decomposition which used classical Poincar\'e inequalities.

\mb
\end{abstract}

\mb {\bf Key-words :}  Calder\'on-Zygmund decomposition, Sobolev spaces, Poincar\'e inequalities.

\tableofcontents

\section{Introduction}

The purpose of this article is to  weaken assumptions of the already known  Calder\'on-Zygmund decomposition for Sobolev functions. This well-known tool was first stated by P. Auscher in \cite{A1}. It exactly corresponds to the Calder\'on-Zygmund decomposition in a context of Sobolev spaces.

\gb Let us briefly recall the ideas of such decomposition. In \cite{Se}, E. Stein stated this decomposition for Lebesgue spaces as following. Let $(X,d,\mu)$ be a space of homogeneous type and $p\geq 1$. Given a function $f\in L^p(X)$, the decomposition gives a precise way of partitioning $X$ into two subsets: one where $f$ is essentially small (bounded in $L^\infty$ norm); the other a countable collection of cubes where $f$ is essentially large, but where some control of the function is obtained in $L^1$ norm. This leads to the associated Calder\'on-Zygmund decomposition of $f$, where $f$ is written as the sum of ``good'' and  ``bad'' functions, using the above subsets.

\mb This decomposition is a basic tool in Harmonic analysis and the study of singular integrals. One of the applications is the following : an $L^2$-bounded Calder\'on-Zygmund operator is of weak type $(1,1)$ and so $L^p$ bounded for every  $p\in(1,\infty)$.

\gb In \cite{A1},  P. Auscher extended these ideas for Sobolev spaces. His decomposition is the following~:

\begin{thm} \label{thm11} Let $n\geq 1$, $p\in[1,\infty)$ and $f\in \mathcal{D}'(\R^n)$ be such that 
$\|\nabla f\|_{L^p} <\infty$. Let $\alpha>0$. Then, one can find a collection of cubes $(Q_i)_i$, functions $g$ and $b_i$  such that 
$$ f= g+\sum_i b_i $$
and the following properties hold:
$$ \|\nabla g\|_{L^\infty} \leq C\alpha, $$
$$ b_i \in W_0^{1,p}(Q_i)\ \text{and} \ \int_{Q_i} |\nabla b_i|^p \le C\alpha^p |Q_i|, $$
$$ \sum_i |Q_i| \le C\alpha^{-p} \int_{\R^n} |\nabla f|^p , $$
$$ \sum_i {\bf 1}_{Q_i} \leq N, $$
where $C$ and  $N$ depend only on the dimension $n$ and on $p$.
\end{thm}

\mb The important point in this decomposition is the fact that the functions $b_i$ are supported in the corresponding balls, while the original Calder\'on-Zygmund decomposition applied to $\nabla f$ would not give this.

\mb The proof relies on an appropriate use of Poincar\'e inequality and was then extended to a doubling manifold with Poincar\'e inequality by P. Auscher and T. Coulhon in \cite{AC}.

\mb This decomposition is used in many works and it appears in various forms and extensions. For example in  \cite{AC} (same proof on manifolds), 
\cite{AM} (on $ \mathbb{R}^n$ but with a doubling weight),  B. Ben Ali's PhD thesis  \cite{Besma} and 
\cite{AB}, \cite{Nadine4} (the Sobolev space is modified to adapt to Schr\"odinger operators), N. 
Badr's PhD thesis \cite{Nadine5} and \cite{Nadine, Nadine1} (used toward interpolation of Sobolev 
spaces on manifolds and measured metric spaces) and in \cite{Nadine3} (Sobolev spaces 
on graphs).

\mb The aim of this article is to extend the proof using other kind of  ``Poincar\'e inequalities''. This work can be integrated in several recent works, where the authors look for replacing the mean-value operators by other ones in the definition of Hardy spaces for example or in the definition of maximal operators (see \cite{BJ,B2,DY,HM,Martell} ... ). Mainly, Section 3 is devoted to the proof of Calder\'on-Zygmund decompositions for Sobolev functions (as in Theorem \ref{thm11}) in an abstract framework of a doubling Riemannian manifold under assumptions involving new kind of Poincar\'e inequalities. Then we give an application to the real interpolation of Sobolev spaces $W^{1,p}$. In Section 4, we focus on a particular case (using the heat semigroup) corresponding to the so-called pseudo-Poincar\'e inequalities. We specify that these new Poincar\'e inequalities are weaker than the classical ones and permit to insure the Calder\'on-Zygmund decomposition for Sobolev functions. We give some applications using this improvement.

\section{Preliminaries}

Throughout this paper we will denote by $\ind_{E}$ the characteristic function of
 a set $E$ and $E^{c}$ the complement of $E$. If $X$ is a metric space, $\Lip$ will be  the set of real Lipschitz functions on $X$ and $\Lip_{0}$ the set of real, compactly supported Lipschitz functions on $X$. For a ball $Q$ in a metric space, $\lambda Q$  denotes the ball co-centered with $Q$ and with radius $\lambda$ times that of $Q$. Finally, $C$ will be a constant that may change from an inequality to another and we will use $u\lesssim
v$ to say that there exists a constant $C$  such that $u\leq Cv$ and $u\simeq v$ to say that $u\lesssim v$ and $v\lesssim u$.

\mb In all this paper, $M$ denotes a complete Riemannian manifold. We write $\mu$ for the Riemannian measure on $M$, $\nabla$ for the
Riemannian gradient, $|\cdot|$ for the length on the tangent space (forgetting the subscript $x$ for simplicity) and
$\|\cdot\|_{L^p}$ for the norm on $ L^p:=L^{p}(M,\mu)$, $1 \leq p\leq +\infty.$  We denote by $Q(x, r)$ the open ball of
center $x\in M $ and radius $r>0$. \\
We will use the positive Laplace-Beltrami operator $\Delta$ defined by
$$ \forall f,g\in C^\infty_0(M), \qquad \langle \Delta f,g\rangle = \langle \nabla f ,\nabla g \rangle.$$
We deal with the Sobolev spaces of order $1$ $W^{1,p}:=W^{1,p}(M)$, where the norm is defined by:
$$ \| f\|_{W^{1,p}(M)} : = \| f\|_p+\|\,||\nabla f|\,\|_{L^p}$$.

\subsection{The doubling property}

\begin{df}[Doubling property] Let $M$ be a Riemannian manifold. One says that $M$ satisfies the doubling property $(D)$ if there exists a constant $C>0$, such that for all $x\in M,\, r>0 $ we have
\begin{equation*}\tag{$D$}
\mu(Q(x,2r))\leq C \mu(Q(x,r)).
\end{equation*}
\end{df}

\begin{lem} Let $M$ be a Riemannian manifold satisfying $(D)$ and let $d=log_{2}C$. Then for all $x,\,y\in M$ and $\theta\geq 1$
\begin{equation}\label{teta}
\mu(Q(x,\theta R))\leq C\theta^{d}\mu(Q(x,R))
\end{equation}
\end{lem} 
\noindent Observe that if $M$ satisfies $(D)$ then
$$ \diam(M)<\infty\Leftrightarrow\,\mu(M)<\infty\,\textrm{ (see \cite{ambrosio1})}. $$
Therefore if $M$ is a complete Riemannian manifold satisfying $(D)$ then $\mu(M)=\infty$.

\begin{thm}[Maximal theorem]\label{MIT} (\cite{coifman2})
Let $M$ be a Riemannian manifold satisfying $(D)$. Denote by $\M$ the uncentered Hardy-Littlewood maximal function
over open balls of $M$ defined by
 $$ \M f(x):=\underset{\genfrac{}{}{0pt}{}{Q \ \textrm{ball}}{x\in Q}} {\sup} \ |f|_{Q} $$
 where $\displaystyle f_{E}:=\aver{E}f d\mu:=\frac{1}{\mu(E)}\int_{E}f d\mu.$
Then for every  $p\in(1,\infty]$, $\M$ is $L^p$ bounded and moreover of weak type $(1,1)$\footnote{ An operator $T$ is of weak type $(p,p)$ if there is $C>0$ such that for any $\alpha>0$, $\mu(\{x;\,|Tf(x)|>\alpha\})\leq \frac{C}{\alpha^p}\|f\|_p^p$.}.
\\
Consequently for $s\in(0,\infty)$, the operator $\M_{s}$ defined by
$$ \M_{s}f(x):=\left[\M(|f|^s)(x) \right]^{1/s} $$
is of weak type $(s,s)$ and $L^p$ bounded for all $p\in(s,\infty]$.
\end{thm}

\subsection{Classical Poincar\'e inequality}
\begin{df}[ Classical Poincar\'{e} inequality on $M$] \label{classP} We say that a complete Riemannian manifold $M$ admits \textbf{a Poincar\'{e} inequality $(P_{q})$} for some $q\in[1,\infty)$ if there exists a constant $C>0$ such that, for every function $f\in \Lip_{0}(M)$\footnote{compaclty supported Lipshitz function defined on $M$.} and every ball $Q$ of $M$ of radius $r>0$, we have
\begin{equation*}\tag{$P_{q}$}
\left(\aver{Q}|f-f_{Q}|^{q} d\mu\right)^{1/q} \leq C r \left(\aver{Q}|\nabla f|^{q}d\mu\right)^{1/q}.
\end{equation*}
\end{df}
\begin{rem} By density of $C_{0}^{\infty}(M)$ in $Lip_0(M)$, we can replace $\Lip_{0}(M)$ by $C_{0}^{\infty}(M)$.
\end{rem}
Let us recall some known facts about Poincar\'{e} inequalities with varying $q$.
 \\
It is known that $(P_{q})$ implies $(P_{p})$ when $p\geq q$ (see \cite{hajlasz4}). Thus, if the set of $q$ such that
$(P_{q})$ holds is not empty, then it is an interval unbounded on the right. A recent result of S. Keith and X. Zhong
(see \cite{KZ}) asserts that this interval is open in $[1,+\infty[$~:

\begin{thm}\label{kz} Let $(X,d,\mu)$ be a complete metric-measure space with $\mu$ doubling
and admitting a Poincar\'{e} inequality $(P_{q})$, for  some $1< q<\infty$.
Then there exists $\epsilon >0$ such that $(X,d,\mu)$ admits
$(P_{p})$ for every $p>q-\epsilon$.
\end{thm}

\subsection{Estimates for the heat kernel} \label{subsec:rapel}


 We recall criterions which rely Poincar\'e inequalities and off-diagonal decays of the heat semigroup. We refer the reader to the work of P. Auscher, T. Coulhon, X. T. Duong and S. Hofmann \cite{ACDH} and \cite{AC} for more details about all these notions and how they are related.
 Let us consider the following two inequalities:
 \begin{equation}\tag{$nhR_p$}
 \|\,\nabla f\,\|_p\leq C (\|\Delta^{\frac{1}{2}}f\|_p+\|f\|_p)
 \end{equation}
 and
  \begin{equation}\tag{$nhRR_p$}
  (\|\Delta^{\frac{1}{2}}f\|_p+\|f\|_p)\leq C \|\,\nabla f\,\|_p.
 \end{equation}
 
\begin{thm} \label{rapel} Let $M$ be a complete doubling Riemannian manifold.
 \begin{itemize}
 \item The inequalities $(nhR_2)$ and $(nhRR_2)$ are always satisfied. 
 \item (\cite{CD1}) Assume that the heat kernel $p_t$ of the semigroup $e^{-t\Delta}$ satisfies the following pointwise estimate:
 \begin{equation}\label{due}
 p_t(x,x) \lesssim \frac{1}{\mu(B(x,t^{1/2}))} \tag{$DUE$}.
 \end{equation}
 Then for all $p\in(1,2]$, $(nhR_p)$ and $(nhRR_{p'})$ hold \footnote{The assumptions in \cite{CD1} are even weaker.}.
 \item (\cite{gri}, Theorem 1.1) Under $(D)$,  $(DUE)$ self-improves into the following Gaussian upper-bound estimate of $p_t$
 \begin{equation}\label{ue}
 p_t(x,y) \lesssim \frac{1}{\mu(B(y,t^{1/2}))}e^{-c\frac{d^{2}(x,y)}{t}} \tag{$UE$}.
\end{equation}
 Note that $(UE)$ implies $(L^1-L^\infty)$ ``off-diagonal'' decays for $(e^{-t\Delta})_{t>0}$.
  \item Under $(UE)$, the collection $(\sqrt{t}\nabla e^{-t\Delta})_{t>0}$ satisfies  ``$L^2-L^2$ off-diagonal decays''.
  \item Under $(DUE)$ and by the analiticity of the heat semigroup, the following pointwise upper bound for the kernel of $\Delta e^{-t\Delta}$: $t \frac{\partial}{\partial t} p_t$ holds (see \cite{Davies2}, Theorem 4 and \cite{gri}, Corollary 3.3):
 \begin{equation} \label{utp} 
 t \left|\frac{\partial }{\partial t} p_t (x,y)\right| \lesssim \frac{1}{\mu(B(y,t^{1/2}))}e^{-c\frac{d^{2}(x,y)}{t}}. 
  \end{equation} 
\end{itemize}
 \end{thm}
 
 \begin{thm}[\cite{LY,saloff}] The conjunction of $(D)$ and Poincar\'e inequality $(P_2)$ on $M$ is equivalent to the following Li-Yau inequality
\begin{equation}\label{ly}
\frac{1}{\mu(B(y,t^{1/2}))}e^{-c_1\frac{d^{2}(x,y)}{t}}\lesssim p_t(x,y) \lesssim \frac{1}{\mu(B(y,t^{1/2}))}e^{-c_2\frac{d^{2}(x,y)}{t}} \tag{$LY$},
\end{equation}
with some constants $c_1,c_2>0$.
\end{thm}
 
\begin{thm}[\cite{ACDH}] The $L^p$-boundedness of the Riesz transform $\nabla (\Delta)^{-1/2}$ implies
\begin{equation} \label{gp}
 \left\|\,|\nabla e^{-t\Delta}|\, \right\|_{L^{p} \to L^{p}} \lesssim \frac{1}{\sqrt{t}} \tag{$G_{p}$}.
 \end{equation}
Moreover, under $(P_2)$ and $(G_{p_0})$  with $p_0>2$,  the collection $(\sqrt{t}\nabla e^{-t\Delta})_{t>0}$ satisfies some  $(L^p-L^p)$ ``off-diagonal'' decays for every $p\in[2,p_0)$.
 \end{thm}

\begin{rem} All these results are proved in their homogeneous version, with homogeneous properties $(R_p)$ and $(RR_p)$.  It is essentially based on the well-known Calder\'on-Zygmund decomposition for Sobolev functions. This tool was extended for non-homogeneous Sobolev spaces (see \cite{Nadine}). Thus by exactly the same proof, we can obtain an analogous
non-homogeneous version  and then prove all these results.
\end{rem}

\subsection{The $K$-method of real interpolation}
\label{subsec:rappel}

We refer the reader  to \cite{bennett}, \cite{bergh} for details on the development of this theory. Here we only recall the essentials to be used in the sequel.

\mb Let $A_{0}$, $A_{1}$ be  two normed vector spaces embedded in a topological Hausdorff vector space $V$. For each  $a\in A_{0}+A_{1}$ and $t>0$, we define the $K$-functional of interpolation by
$$
K(a,t,A_{0},A_{1})=\displaystyle \inf_{a
=a_{0}+a_{1}}(\| a_{0}\|_{A_{0}}+t\|
a_{1}\|_{A_{1}}).
$$

\gb For $0<\theta< 1$, $1\leq q\leq \infty$, we denote by $(A_{0},A_{1})_{\theta,q}$ the interpolation space between $A_{0}$ and $A_{1}$:
\begin{displaymath}
    (A_{0},A_{1})_{\theta,q}=\left\lbrace a \in A_{0}+A_{1}:\|a\|_{\theta,q}=\left(\int_{0}^{\infty}(t^{-\theta}K(a,t,A_{0},A_{1}))^{q}\,\frac{dt}{t}\right)^{\frac{1}{q}}<\infty\right\rbrace.
\end{displaymath}
It is an exact interpolation space of exponent $\theta$ between $A_{0}$ and $A_{1}$ (see \cite{bergh}, Chapter II).
\begin{df}
Let $f$  be a measurable function on a measure space $(X,\mu)$. The decreasing rearrangement of $f$ is the function $f^{*}$ defined for every $t\geq 0$ by
$$
f^{*}(t)=\inf \left\lbrace\lambda :\, \mu (\left\lbrace x:\,|f(x)|>\lambda\right\rbrace)\leq
t\right\rbrace.
$$
The maximal decreasing rearrangement of
$f$ is the function $f^{**}$ defined for every $t>0$ by
$$
f^{**}(t)=\frac{1}{t}\int_{0}^{t}f^{*}(s) ds.
$$
\end{df}

\begin{prop} \label{prop}
>From the properties of $f^{**}$, we mention:
\begin{itemize}
\item[1.] $(f+g)^{**}\leq f^{**}+g^{**}$.
\item[2.] $(\mathcal{M}f)^{*}\sim f^{**}$.
\item[3.] $\mu(\left\lbrace x;\,|f(x)|>f^{*}(t)\right\rbrace )\leq t$.
\item[4.] $\forall p\in(1,\infty]$, $\|f^{**}\|_{p}\sim\|f\|_{p}$.
\end{itemize}
\end{prop}

\mb We exactly know the functional $K$ for Lebesgue spaces~:
\begin{prop} \label{prop2} Take $0<p_0<p_1\leq \infty$. We have~:
$$K(f,t,L^{p_0},L^{p_1}) \simeq \left(\int_0^{t^{\alpha}} \left[f^{*}(s)\right]^{p_0} ds \right)^{1/p_0} + t \left(\int_{t^{\alpha}}^\infty \left[f^{*}(s)\right]^{p_1} ds \right)^{1/p_1},$$
where $\frac{1}{\alpha}=\frac{1}{p_0}-\frac{1}{p_1}$.
\end{prop}

\section{New ``Calder\'on-Zygmund'' decompositions for Sobolev functions.}
\label{sec:methode2}

In the introduction, we recalled the main use of ``Calder\'on-Zygmund'' decompositions for Sobolev functions. In the previously cited works, this decomposition relies on Poincar\'e inequalities and some ``tricks'' with the mean-value operators. We present here similar arguments with abstract operators, requiring new ``Poincar\'e inequalities''. Then, we  give some applications to real interpolation of Sobolev spaces.

\subsection{Decomposition using abstract ``oscillation operators''}

Let $\A:=(A_Q)_Q$ be a collection of operators (acting from $W^{1,p}$ to $W^{1,p}_{loc}$) indexed by the balls of the manifold ($A_Q$ can be thought to be similar to the mean operator over the ball $Q$).

\begin{df} \label{def:opeM} We define a new maximal operator associated to this collection: for $1\leq s \leq p \leq \infty$ and all functions $f\in  W^{1,p}$
$$M_{\A,s}(f)(x) := \sup_{Q;\,Q\ni x} \ \frac{1}{\mu(Q)^{1/s}} \left\| A_Q(f)\right\|_{W^{1,s}(Q)}.$$
\end{df}

\mb Let us now define the  assumptions that we need on the collection $\A$.

\begin{df} $1)$ We say that for $q\in[1,\infty]$ \footnote{we take the supremun instead of the $L^q$ average when $q=\infty$.}, the manifold $M$ satisfies a Poincar\'e inequality $(P_q)$ relatively to the collection $\A$ if there is a constant $C$ such that for every ball $Q$ (of radius $r_Q$) and for all functions $f\in W^{1,p};\, p\geq q$:
\begin{align*}
\left(\aver{Q} \left|f-A_Q(f)\right|^q d\mu \right)^{1/q} \leq C r_Q \sup_{s\geq 1}
\left(\aver{sQ}\left(|f|+|\nabla f|\right)^q d\mu \right)^{1/q}.
\end{align*}

\mb $2)$ For $1\leq q \leq r\leq \infty$, we say that the collection  $\A$ satisfies ``$L^q-L^r$ off-diagonal estimates'' if
 \begin{itemize}
 \item[a.] there are constants $C'>0$ and $N\in {\mathbb N}^*$ such that for all equivalent balls $Q,\,Q'$ (i.e. $Q\subset Q' \subset N Q$) and all functions $f\in W^{1,p};\, p\geq q$, we have
\be{off2} \frac{1}{\mu(Q)^{1/r}} \left\|A_{Q}(f)-A_{Q'}(f) \right\|_{L^r(NQ)} \leq C' r_{Q} \inf_{NQ} \M_{q}\left(|f|+|\nabla f|\right) \ee
  \item[b.] and  for every ball $Q$
  \begin{align}
\frac{1}{\mu(Q)^{1/r}}\left\|A_{Q}(f)\right\|_{W^{1,r}(Q)} \leq C' \inf_{Q} \M_{q}\left(|f|+|\nabla f|\right). \label{off}
\end{align}
\end{itemize}
\end{df}

\mb Here is our main result~:

\begin{thm}\label{thm:CZ} Let $M$ be a complete Riemannian manifold satisfying  $(D)$ and of infinite measure. Consider a collection $\A =(A_Q)_Q$ of operators defined on $M$. Assume that $M$ satisfies the Poincar\'e inequality $(P_{q})$  relatively to the collection $\A$ for some $q\in[1,\infty)$,  and that $\A$ satisfies ``$L^q-L^r$ off-diagonal estimates'' for some $r\in(q,\infty]$. \\
 Let $q\leq p<r$, $f\in W^{1,p}$ and $\alpha>0$. Then one can find a collection of balls $(Q_{i})$, functions $g\in W^{1,r}$ and $b_{i}\in W^{1,q}$ with the following properties
\begin{equation}
f = g+\sum_{i}b_{i} \label{df}
\end{equation}
\begin{equation}
\left\|g\right\|_{W^{1,r}} \lesssim \|f\|_{W^{1,p}}^{p/r} \alpha^{1-p/r}, \; \int_{\cup_i{Q_i}}(|g|^{r}+|\nabla g|^{r})d\mu\lesssim \alpha^{r}\mu(\cup_iQ_i) \label{eg}
\end{equation}
\begin{equation}
supp \,(b_{i})\subset Q_{i}, \,\left\|b_{i}\right\|_{W^{1,q}} \lesssim \alpha \mu(Q_{i})^{1/q} \label{eb}
\end{equation}
\begin{equation}
\sum_{i}\mu(Q_{i})\leq C\alpha^{-p}\int (|f|+|\nabla f|)^{p} d\mu
\label{eB}
\end{equation}
\begin{equation}
\sum_{i}{\bf 1}_{Q_{i}}\leq N \label{rb}.
\end{equation}
 \end{thm}

\begin{rem} From the assumed ``$L^q-L^r$ off-diagonal estimates'' for $\A$ and Theorem \ref{MIT}, we deduce that the maximal operator $M_{\A,q}$ is continuous from $W^{1,q}$ to $L^{q,\infty}$ and from $W^{1,p}$ to $L^{p}$ for $p\in (q,r]$. 
\end{rem}

 \dem We follow the ideas of \cite{Nadine} where the result is proved for the particular case
 $$A_Q(f):= \aver {Q} f d \mu.$$ \\
Let $f\in W^{1,p} $ and $\alpha>0$. Consider the set
$$ \Omega :=\left\{ x\in M;\ \M_{q}(|f|+|\nabla f|)(x)+M_{\A,q}(f)(x) >\alpha\right\}.$$
We can assume that this set is non empty (otherwise the result is obvious taking $g=f$). With this assumption, the different maximal operators are of ``weak type $(p,p)$'' so
 \begin{align}
    \mu(\Omega)& \leq C \alpha^{-p} \left(\int | f|^{p} d\mu +\int |\nabla f|^{p} d\mu\right) \label{mO} \\
    & <+\infty. \nonumber
\end{align}
   In particular $\Omega\neq M$ as $\mu(M)=\infty$. Let $F$ be the complement of $\Omega$. Since $\Omega$ is an open set distinct of $M$, we can take $(\underline{Q_{i}})$ a Whitney decomposition of $\Omega$. That is the balls  $\underline{Q_{i}}$ are pairwise disjoint and there exist two constants $C_{2}>C_{1}>1$, depending only
on the metric, such that
\begin{itemize}
\item[1.] $\Omega=\cup_{i}Q_{i}$ with $Q_{i}=
C_{1}\underline{Q_{i}}$ and the balls $Q_{i}$ have the bounded overlap property;
\item[2.] $r_{i}=r(Q_{i})=\frac{1}{2}d(x_{i},F)$ and $x_{i}$ is
the center of $Q_{i}$;
\item[3.] each ball $C_{2}\underline{Q_{i}}$ intersects $F$ ($C_{2}=4C_{1}$ works) and  we define $\overline{Q_{i}}=2C_{2}\underline{Q_{i}}$.
\end{itemize}
For $x\in \Omega$, denote $I_{x}=\left\lbrace i:x\in Q_{i}\right\rbrace$. By the bounded overlap property of the balls $Q_{i}$, we have that $\sharp I_{x} \leq N$ with a numerical integer $N$. Fixing $j\in I_{x}$ and using the properties of the $Q_{i}$'s, we easily see that $\frac{1}{3}r_{i}\leq r_{j}\leq 3r_{i}$ for all $i\in I_{x}$. In particular, $Q_{i}\subset 7Q_{j}$ for all $i\in I_{x}$.

\mb Condition (\ref{rb}) is nothing but the bounded overlap property of the $Q_{i}$'s  and (\ref{eB}) follows from
(\ref{rb}) and  (\ref{mO}).

\gb Observe that the doubling property and the fact that $\overline{Q_{i}} \cap F\neq \emptyset$ yield
\begin{align}
\int_{Q_{i}} (|f|^{q}+|\nabla f|^{q}+\left|A_{\overline{Q_{i}}}f\right|^{q}+ \left|\nabla A_{\overline{Q_{i}}}f\right|^{q}) d\mu & \leq
\int_{\overline{Q_{i}}} (|f|^{q}+|\nabla f|^{q}+\left|A_{\overline{Q_{i}}}f\right|^{q}+ \left|\nabla A_{\overline{Q_{i}}}f\right|^{q}) d\mu \nonumber \\
 & \leq \inf_{\overline{Q_{i}}} \left[\M_{q}(|f|+|\nabla f|) + M_{\A,q}(f)\right]^q \mu(\overline{Q_{i}}) \nonumber \\
 & \leq \alpha^{q} \mu(\overline{Q_{i}}) \nonumber \\
 & \lesssim \alpha^{q}\mu(Q_{i}). \label{f}
\end{align}
We now define the functions $b_{i}$. Let $(\chi_{i})_{i}$ be  a partition of unity of $\Omega$ associated to the
covering $(\underline{Q_{i}})$, such that for all $i$, $\chi_{i}$ is a Lipschitz function supported in $Q_{i}$ with
$\displaystyle\|\,|\nabla \chi_{i}|\, \|_{\infty}\lesssim r_{i}^{-1}$. Set
$$b_{i}:=(f-A_{\overline{Q_i}}f)\chi_{i}.$$
It is clear that $\textrm{supp}(b_{i}) \subset Q_{i}$.
Let us estimate $\left\|b_{i}\right\|_{W^{1,q}(Q_i)}$. We have
\begin{align*}
\int_{Q_{i}} |b_{i}|^{q} d\mu
&=\int_{Q_{i}} \left|(f-A_{\overline{Q_i}}(f)\right|^q d\mu
\\
&\lesssim
\int_{Q_{i}}|f|^{q}d\mu+\int_{Q_{i}}\left|A_{\overline{Q_i}}(f)\right|^{q} d\mu
\\
&\lesssim \alpha^{q} \mu(Q_{i}).
\end{align*}
We applied (\ref{f}) in the last inequality. Since $$\nabla\Bigl((f-A_{\overline{Q_i}}f)\chi_{i}\Bigr)=\chi_{i}\left(\nabla f- \nabla A_{\overline{Q_i}}f \right)+\left(f-A_{\overline{Q_i}}f\right)\nabla\chi_{i},$$
we have
\begin{align*}
\int_{Q_{i}} |\nabla b_{i}|^{q} d\mu
\lesssim \int_{Q_{i}} \left|\nabla f-\nabla A_{\overline{Q_i}}(f)\right|^q d\mu + \frac{1}{r_i^q} \int_{Q_i} \left|f-A_{\overline{Q_i}}f\right|^q d\mu.
\end{align*}
The first term is estimated as above for $b_i$. Thus
$$\int_{Q_{i}} \left|\nabla f-\nabla A_{\overline{Q_i}}(f)\right|^q d\mu \lesssim \alpha^q \mu(Q_i).$$
For the second term, the Poincar\'e inequality $(P_q)$ (relatively to the collection $\A$) shows that
\begin{align*}
\frac{1}{r_i^q}\int_{Q_{i}} \left|f-A_{\overline{Q_i}}(f)\right|^q d\mu &  \lesssim \sup_{s\geq 1}
\frac{\mu(Q_i)}{\mu(sQ_i)} \int_{sQ_{i}} (|f|^q+|\nabla
f|^{q})d\mu \\
&\lesssim \alpha^{q}\mu(Q_{i}).
\end{align*}
We used that for all $s\geq 1$, $s\overline{Q_i}$ meets $F$ and (\ref{f}) for $sQ_i$ instead of $Q_i$.
Therefore (\ref{eb}) is proved.

\gb Set $\displaystyle g=f-\sum_{i}b_{i}$, then it remains to prove (\ref{eg}). Since the sum is locally finite on
$\Omega$, $g$ is defined  almost everywhere on $M$ and $g=f$ on $F$. Observe that $g$ is a locally integrable function
on $M$. This follows from the fact that $b=f-g\in L^{q}$ here (for the homogeneous case, one can easily prove that $b\in L^1_{loc}$). Note that $\displaystyle \sum_{i}\chi_{i}={\bf 1}_{\Omega}$ and
$\displaystyle \sum_{i}\nabla\chi_{i}= \nabla {\bf 1}_{\Omega}$. We then have
\begin{align}
\nabla g & = \nabla f -\sum_{i}\nabla b_{i} \nonumber \\
& = \nabla f-\left(\sum_{i} \chi_{i}\left[\nabla f - \nabla A_{\overline{Q_{i}}}f \right]\right) -\sum_{i}(f-A_{\overline{Q_{i}}}(f))\nabla \chi_{i} \nonumber \\
& ={\bf 1}_{F}(\nabla f) + \sum_{i} \chi_{i}\nabla A_{\overline{Q_{i}}}f -\sum_{i}A_{\overline{Q_{i}}}(f)\nabla\chi_{i}
- f \nabla {\bf 1}_{\Omega}. \label{eqal}
\end{align}
The definition of $F$ and the Lebesgue differentiation theorem yield ${\bf 1}_{F}(|f|+|\nabla f|)\leq \alpha\;\mu -$a.e. We deduce that (with an interpolation inequality) for $\frac{1}{r}=\frac{\theta}{p}$~:
\begin{align*}
 \left\| {\bf 1}_{F}(|f|+|\nabla f|) \right\|_{L^r} & \lesssim \left\| {\bf 1}_{F}(|f|+|\nabla f|) \right\|_{L^p}^{\theta} \left\| {\bf 1}_{F}(|f|+|\nabla f|) \right\|_{L^\infty}^{(1-\theta)} \\
 & \lesssim \|f\|_{W^{1,p}}^{p/r} \alpha^{1-p/r}.
\end{align*}
We control the second term in (\ref{eqal}) using  the ``off-diagonal'' decays of $\A$: (\ref{off}).
 We recall that $\overline{Q_{i}} = 2C_2 \underline{Q_i}$. We deduce that
\begin{align}
\left\|\,|\nabla A_{\overline{Q_{i}}}f|\,\right\|_{L^r(Q_i)} & \lesssim \mu(Q_i)^{1/r} \inf_{\overline{Q_{i}}} \M_{q}\left(|f|+|\nabla f|\right) \nonumber \\
 & \lesssim \alpha \mu(Q_i)^{1/r}. \label{oui}
\end{align}
The last inequality is due to the fact that $\overline{Q_{i}}\cap F\neq \emptyset$. Then the bounded overlap property of the covering $(Q_i)_i$ gives us
\begin{align*}
\left\| \sum_{i} \chi_{i}(x)|\nabla A_{\overline{Q_{i}}}f |\,\right\|_{L^r} & \lesssim \left( \sum_i \left\|\,| \nabla A_{\overline{Q_{i}}}f |\,\right\|_{L^r(Q_i)} ^r \right)^{1/r} \\
 & \lesssim \left( \alpha^r \sum_i \mu(Q_i) \right)^{1/r} \\
  & \lesssim \alpha (\mu(\Omega))^{1/r}
\end{align*}
We claim that a similar estimate holds for $h=\sum_{i}\left[A_{\overline{Q_{i}}}(f)-f\right] \nabla \chi_{i}$~: we have
$\|h \|_{L^r}\lesssim \alpha (\mu(\Omega))^{1/r}$.
 \\ 
To prove this, we fix a point $x\in \Omega$ and let $Q_{j}$  be a Whitney
ball containing $x$. For all $i\in I_{x}$ as $r_{Q_i} \simeq r_{Q_j}$, we have \be{g}
\left\|A_{\overline{Q_{i}}}(f)-A_{\overline{Q_{j}}}(f) \right\|_{L^r(Q_i)} \lesssim r_{j}\mu(Q_j)^{1/r} \alpha.\ee
Indeed, since $Q_{i} \subset 7Q_{j}$, this is a direct consequence of the assumed ``off-diagonal'' decays and the fact that
$10\overline{Q_{i}}\cap F \neq \emptyset$. Using $\displaystyle \sum_{i}\nabla\chi_{i}(x)=0$, we deduce that \be{eq:CZ}
\left\|h \right\|_{L^r(Q_j)} \lesssim \sum_{i\in I_x} \left\|A_{\overline{Q_{i}}}(f)-A_{\overline{Q_{j}}}(f)
\right\|_{L^r(Q_j)} r_j^{-1}
 \lesssim N\alpha \mu(Q_j)^{1/r} \lesssim \alpha\mu(Q_j)^{1/r}. \ee
Using again the bounded overlap property of the $(Q_i)_i$'s, it follows that
$$ \|h \|_{L^r}\lesssim \alpha (\mu(\Omega))^{1/r}.$$
Hence
$$ \| \,|\nabla g|\, \|_{L^r(\Omega)}\lesssim \alpha (\mu(\Omega))^{1/r}.$$
 Then (\ref{eB}) and the $L^r$ estimate of $|\nabla g|$ on $F$ yield $\| \nabla g \|_{L^r}\lesssim \|f\|_{W^{1,p}}^{p/r}\alpha^{1-p/r}$.
Let us now estimate $\| g \|_{L^r}$. We have $\displaystyle g=f{\bf 1}_{F}+\sum_{i}A_{\overline{Q_{i}}}(f)\chi_{i}$. Since $|f|{\bf 1}_{F}\leq \alpha$, still need to estimate $\|\sum_{i}A_{\overline{Q_{i}}}(f)\chi_{i}\|_{L^r}$. Note that as in (\ref{oui}), we similarly have~ for every $i$
 \be{f_B}
 \left\|A_{\overline{Q_{i}}}(f)\right\|_{L^r(Q_i)} \lesssim \alpha \mu(Q_i)^{1/r}.
 \ee
As above, this last inequality yields (thanks to the bounded overlap property of the $(Q_i)_i$)
$$\| g \|_{L^r(\Omega)}\lesssim \alpha (\mu(\Omega))^{1/r}.$$ 
Finally,  (\ref{eB}) and the $L^r$ estimate of $g$ on $F$ yield
 $\| g \|_{L^r}\lesssim \|f\|_{W^{1,p}}^{p/r} \alpha^{1-p/r}$.
Therefore we proved that $g$ belongs to $W^{1,r}$ with the desired boundedness.
\findem

\begin{rem} \label{rem:CZ} Note that in this decomposition, $\nabla {\bf 1}_{\Omega}$
corresponds to a singular distribution, supported in $\partial \Omega$. In the previous proof, we  considered that
the distribution $\nabla {\bf 1}_{\Omega}$ corresponds to a function, vanishing almost everywhere. The estimate
(\ref{eq:CZ}) shows that $h$ (considered as an $L^1_{loc}$-function) satisfies the good property. We also have to check
that $h$ can be considered as an $L^1_{loc}$-function. This is due to the following fact~
$$ \sum_{i,j}\left[A_{\overline{Q_{j}}}(f)\chi_j -f\right] \nabla \chi_{i} =0$$
in the distributional sense. This equality shows that when we are close to  ${\textrm{supp}}(\sum \nabla \chi_i)= \partial \Omega$, the corresponding operator $A_{\overline{Q_{j}}}$ tends to the identity operator, due to Poincar\'e inequality.
We do not detail this technical problem and refer to \cite{A2}.
\end{rem}

\begin{rem}\label{part} In the case where the operator $A_Q$ is the mean-operator over the ball $Q$, the assumption ``$M_{\A,q}=\M_{q}$ is continuous from $W^{1,p}$ to $L^{p,\infty}$'' is always satisfied.
The Poincar\'e inequality $(P_q)$ corresponds to the ``classical one'' (in fact it is weaker since that in the classical one it appears only the $L^q(Q)$ norm of the gradient of the function) .
Moreover ``$L^q-L^\infty$ off-diagonal estimates'' hold obviously. Thus, we regain the well-known Calder\'on-Zygmund decomposition in Sobolev spaces.
\end{rem}

\subsection{Application to real Interpolation of Sobolev spaces.}

\mb As described in \cite{Nadine1}, such a ``Calder\'on-Zygmund'' decomposition in Sobolev spaces is sufficient to
obtain a real interpolation result for Sobolev spaces.

\begin{thm} \label{thm:int} Let $M$ be a complete Riemannian manifold of infinite measure satisfying $(D)$ and admitting a Poincar\'e inequality $(P_{q})$ for some $q\in[1,\infty)$ relatively to the collection $\A$. Assume that $\A$ satisfies ``$L^q-L^r$ off-diagonal estimates'' for an $r\in(q,\infty]$. 
Then for  $1\leq s\leq p<r\leq \infty$ with $p>q$, the space $W^{1,p}$ is a real interpolation space between $W^{1,s}$ and $W^{1,r}$. More precisely
$$ W^{1,p}=(W^{1,s},W^{1,r})_{\theta,p}$$
where $\theta\in(0,1)$ such that
$$ \frac{1}{p}:=\frac{1-\theta}{s}+\frac{\theta}{r}<\frac{1}{q}.$$
\end{thm}

\mb We do not detail the proof and refer the reader to \cite{Nadine1} for the link between  such a ``Calder\'on-Zygmund'' decomposition and interpolation results. We briefly explain the main steps of the proof.

\dem It is sufficient to prove that there exists $C>0$ such that for every $f\in W^{1,p}$ and $t>0$,
\begin{align} 
\lefteqn{K(f,t, W^{1,s},W^{1,r})} & & \nonumber \\
 & & \lesssim\left( t^{\frac{r}{r-s}} \left[|f|^{q**}+|\nabla f|^{q**}\right]^{1/q}(t^{\frac{rs}{r-s}})+t\left[\int_{t^{\frac{rs}{r-s}}}^\infty \left(\M (|f|+|\nabla f|)^{q}\right)^{*r/q}(u)du \right]^{1/r} \right). \label{eq:am}
\end{align}
We consider the previous Calder\'{o}n-Zygmund decomposition for $f$ with
$$\alpha=\alpha(t)=\left[\M_q (|f|+|\nabla f|) + M_{A,q}(f)\right]^{q*\frac{1}{q}}(t^{\frac{rs}{r-s}}).$$
We write $ \displaystyle f=\sum_{i}b_{i}+g=b+g $ where
$(b_{i})_{i},\,g$ satisfy the properties of Theorem \ref{thm:CZ}. From the bounded overlap property of the $B_{i}$'s, it follows that 
\begin{align*}
\| b \|_{W^{1,s}}^s& \leq N\sum_{i}
\|b_{i}\|_{W^{1,s}}^s
\\
&\lesssim \alpha^{s}(t)\sum_{i}\mu(B_{i})
\\
&\lesssim \alpha^{s}(t)\mu(\Omega_t),
\end{align*}
with $\Omega_t=\cup_i B_i$.
For $g$, we have as in \cite{Nadine1}, proof of Theorem 4.2, p.15 
\begin{align*}
\int_{F_t}(|g|^r+|\nabla g|^{r})\,d\mu&=\int_{F_t}(|f|^r+|\nabla f|^{r})\,d\mu
\\
&\lesssim \int_{t^{\frac{rs}{r-s}}}^{\infty}(\mathcal{M}(|f|+|\nabla f|)^{q})^{*\frac{r}{q}}(u)du
\\
&+t^{\frac{rs}{r-s}}(|f|^{q**}+|\nabla f|^{q**})^{\frac{r}{q}}(t^{\frac{rs}{r-s}})
\end{align*}
where $F_t$ is the complement of $\Omega_t$. For the Sobolev norm of $g$ in $\Omega$, we use the estimate of the Calder\'on-Zygmund decomposition. Moreover, since $(\M f)^{*}\sim f^{**}$ and $(f+g)^{**}\leq f^{**}+g^{**}$ (c.f \cite{bennett},\cite{bergh}) and thanks to the ``$(L^q-L^r)$ off-diagonal'' assumption on $\A$, we have $$
\alpha(t)\lesssim \left(|f|^{q**{\frac{1}{q}}}(t^{\frac{rs}{r-s}})+|\nabla
f|^{q**{\frac{1}{q}}}(t^{\frac{rs}{r-s}})\right).$$
The choice of $\alpha(t)$ implies $\mu(\Omega_t)\leq t^{\frac{rs}{r-s}}$ (c.f \cite{bennett},\cite{bergh}). 
Finally (\ref{eq:am}) follows from the fact that
$$ K(f,t,W^{1,s},W^{1,r}) \leq \|b\|_{W^{1,s}} + t\|g\|_{W^{1,r}}$$
and the good estimates of $\|b\|_{W^{1,s}}$ and $\|g\|_{W^{1,r}}$.
\findem

\begin{rem} \label{rem:local} As explained in \cite{Nadine,Nadine1}, to interpolate the non-homogeneous Sobolev spaces, it is sufficient  to assume local doubling $(D_{loc})$ and local Poincar\'e inequality $(P_{qloc})$ relatively to $\A$. In these assumptions, we restrict to balls $Q$ of radius sufficiently small.
\end{rem}

\mb We now give an homogeneous version of all these results and then give applications.
 
\subsection{Homogeneous version}
We begin recalling the definition of homogeneous Sobolev spaces on a manifold.

\gb Let $M$ be a $C^{\infty}$ Riemannian manifold of dimension $n$.
For $1\leq p\leq \infty$, we define $\overset{.}{E^{1,p}}$ to be the vector space of distributions $ \varphi $
 with $|\nabla \varphi |\in L^{p}$, where $\nabla \varphi$ is the distributional gradient of $ \varphi$. We equip $\overset{.}{E^{1,p}}$  with the semi-norm 
 $$ 
 \|\varphi\|_{\overset{.}{E^{1,p}}}=\|\,|\nabla \varphi|\,\|_{L^p}.
 $$
The homogeneous Sobolev space $\overset{.}{W^{1,p}}$ is then the quotient space $\overset{.}{E_{p}^{1}}/\mathbb{R}$.

\begin{rem} 1.For all $\varphi \in \overset{.}{E^{1,p}}$, \, $\|\overline{\varphi}\|_{\overset{.}{W^{1,p}}}=\|\,|\nabla \varphi|\,\|_{L^p}$, where $\overline{\varphi}$ denotes the class of $\varphi$.
\\
2. The space $\overset{.}{W^{1,p}}$ is a Banach space (see \cite{goldshtein}).
\end{rem}

\mb We then have all the homogeneous version of our results. We only state them, their proofs being the same as in the non-homogeneous case with few modifications due to the homogeneous norm.

\mb Let $\A:=(A_Q)_Q$ be a collection of operators (acting from $\dot{W}^{1,p}$ to $\dot{W}^{1,p}_{loc}$) indexed by the balls of the manifold.
 We define analogously new homogeneous maximal operator associated to this collection: for $1\leq s \leq p \leq \infty$ and all functions $f\in  \dot{W}^{1,p}$
$$\dot{M}_{\A,s}(f)(x) := \sup_{Q;\,Q\ni x} \ \frac{1}{\mu(Q)^{1/s}} \left\| \,| \nabla A_Q(f)|\,\right\|_{L^s(Q)}.$$
\\
The  assumptions that we need on the collection $\A$ are then the following:

\begin{df} $1)$ We say that for $q\in[1,\infty]$, the manifold $M$ satisfies an homogeneous Poincar\'e inequality $(\dot{P}_q)$ relatively to the collection $\A$ if there is a constant $C$ such that for every ball $Q$ (of radius $r_Q$) and for all functions $f\in \dot{W}^{1,p};\, p\geq q$:
\begin{align*}
\left(\aver{Q }\left|f-A_Q(f)\right|^q d\mu \right)^{1/q} \leq C r_Q \sup_{s\geq 1}
\left(\aver{sQ}|\nabla f|^q d\mu \right)^{1/q}.
\end{align*}

\mb $2)$ We say that the collection  $\A$ satisfies ``$L^q-L^r$ homogeneous off-diagonal estimates'' if
 \begin{itemize}
 \item[a.] there are constants $C'>0$ and $N\in {\mathbb N}^*$ such that for all equivalent balls $Q,\,Q'$ (i.e. $Q\subset Q' \subset N Q$; $N\in\mathbb{N}^*$) and all functions $f\in \dot{W}^{1,p};\, p\geq q$, we have
$$ \frac{1}{\mu(Q)^{1/r}} \left\|A_{Q}(f)-A_{Q'}(f) \right\|_{L^r(NQ)} \leq C' r_{Q} \inf_{NQ} \M_{q}\left(|\nabla f|\right)$$
  \item[b.] and  for every ball $Q$
  \begin{align}
\frac{1}{\mu(Q)^{1/r}}\left\| \,| \nabla A_{Q}(f)|\,\right\|_{L^r( Q)} \leq C' \inf_{Q} \M_{q}\left(|\nabla f|\right). \label{offh}
\end{align}
\end{itemize}
\end{df}

\mb Then, we get the homogeneous version of the Calder\'on-Zygmund decomposition:

\begin{thm}\label{thm:CZH} Let $M$ be a complete Riemannian manifold satisfying  $(D)$ and of infinite measure. Consider a collection $\A =(A_Q)_Q$ of operators defined on $M$.
 Assume that $M$ satisfies the Poincar\'e inequality $(\dot{P}_{q})$  relatively to the collection $\A$ for some $q\in[1,\infty)$  and that $\A$ satisfies $L^q-L^r$ `` homogeneous off-diagonal estimates'' for an $r\in(q,\infty]$. \\ 
 Let $f\in \dot{W}^{1,p}$ and $\alpha>0$. Then one can find a collection of balls $(Q_{i})$, functions $g\in \dot{W}^{1,r}$ and $b_{i}\in \dot{W}^{1,q}$ with the following properties
\begin{equation}
f = g+\sum_{i}b_{i} \label{dfh}
\end{equation}
\begin{equation}
\left\|g\right\|_{\dot{W}^{1,r}} \lesssim \|f\|_{\dot{W}^{1,p}}^{p/r} \alpha^{1-p/r}, \; \int_{\cup_i{Q_i}}|\nabla g|^{r}d\mu\lesssim \alpha^{r}\mu(\cup_iQ_i) \label{egh}
\end{equation}
\begin{equation}
supp \,(b_{i})\subset Q_{i}, \,\left\|b_{i}\right\|_{\dot{W}^{1,q}} \lesssim \alpha \mu(Q_{i})^{1/q} \label{ebh}
\end{equation}
\begin{equation}
\sum_{i}\mu(Q_{i})\leq C\alpha^{-p}\int |\nabla f|^{p} d\mu
\label{eBh}
\end{equation}
\begin{equation}
\sum_{i}{\bf 1}_{Q_{i}}\leq N \label{rbh}.
\end{equation}
 \end{thm}

\mb This decomposition will give us the following homogeneous interpolation result:

\begin{thm} \label{ISH} Let $M$ be a complete Riemannian manifold of infinite measure satisfying  $(D)$ and admitting a Poincar\'e inequality $(\dot{P}_{q})$ for some $q\in[1,\infty)$ relatively to the collection $\A$.  Assume that $\A$ satisfies $L^q-L^r$ `` homogeneous off-diagonal estimates'' for an $r\in(q,\infty]$.  \\
Then for  $1\leq s \leq p<r\leq \infty$ with $p>q$, 
the space $\dot{W}^{1,p}$ is a real interpolation space between $\dot{W}^{1,s}$ and $\dot{W}^{1,r}$. More precisely
$$ \dot{W}^{1,p}=(\dot{W}^{1,s},\dot{W}^{1,r})_{\theta,p}$$
where $\theta\in(0,1)$ such that
$$ \frac{1}{p}:=\frac{1-\theta}{s}+\frac{\theta}{r}<\frac{1}{q}.$$
\end{thm}
\section{Pseudo-Poincar\'e inequalities and Applications}

\subsection{The particular case of ``Pseudo-Poincar\'e Inequalities''}

Thanks to \cite{A1,A}, we know that under $(D)$, a Poincar\'e inequality $(P_q)$ guarantees the assumptions of Theorem \ref{thm:CZ} when $A_Q$ is the  mean-operator over the ball $Q$. Thus it permits to prove a Calder\'on-Zygmund decomposition for Sobolev functions. \\
The aim of this subsection is to show, using a particular choice of operators $A_Q$, that our assumptions are weaker than the classical Poincar\'e inequality used in the already known decomposition.  \\
Let $\Delta$ be the positive Laplace-Beltrami operator and let us set $A_Q:=e^{-r_Q^2\Delta}$ for each ball $Q$ of radius $r_Q$. In all this section, we work with these operators. In order to obtain a Calder\'on-Zygmund decomposition as in Theorem \ref{thm:CZ}, we need to put some assumptions on $(A_Q)_Q$ as those in Section \ref{sec:methode2}.

\mb According to this choice of operators, we define what are ``Pseudo-Poincar\'e inequalities''.

\begin{df}[Pseudo-Poincar\'{e} inequality on $M$] We say that a complete Riemannian manifold $M$ admits \textbf{a pseudo-Poincar\'{e} inequality $(\widetilde{P}_{q})$}
for some $q\in[1,\infty)$ if there exists a constant $C>0$ such that, for every function $f\in C^\infty_0$ and every ball
$Q$ of $M$ of radius $r>0$, we have
\begin{equation*} \tag{$\tilde{P}_{q}$}
\left(\aver{Q}|f-e^{-r^2\Delta}f|^{q} d\mu\right)^{1/q} \leq C r \sup_{s\geq 1} \left(\aver{sQ} |\nabla
f|^{q}d\mu\right)^{1/q}. 
\end{equation*}
\end{df}

\mb Pseudo-Poincar\'e inequalities  corresponds to what we called  Poincar\'e inequality relatively to this collection $\A$ (the homogeneous version, we can also consider the non-homogeneous one).

\mb We begin showing that pseudo-Poincar\'{e} inequalities are implied by the classical Poincar\'e inequalities. We denote  
\be{q0} q_0:= \inf \{ q\in [1,\infty); \ (P_q) \;\textrm{ holds } \}. \tag{$q_0$} \ee 

\begin{prop} \label{prop:pseudo} Let $M$ be a complete manifold satisfying $(D)$ and admitting a Poincar\'e inequality $(P_{q})$ for some $1\leq q<\infty$. 
\begin{itemize}
\item[1.] If $q_0<2$ then the pseudo-Poincar\'e inequality ($\widetilde{P}_{q}$) holds.
\item[2.] If $q_0\geq 2$, we moreover assume $(DUE)$. Then ($\widetilde{P}_{q}$) also holds.
\end{itemize}
\end{prop}

\mb Before proving this proposition, we give the following covering Lemma.
\begin{lem} \label{cov} Let  $M$ be a complete manifold satisfying $(D)$. Let $Q$  a ball of radius $r_Q$. Then there exists a bounded covering $(Q_j)_j$ of $Q$ with balls of radius $t^{1/2}$ for $0<t\leq r_Q^2$.  Moreover,  for $s\geq 1$, the collection $(sQ_j)_j$ is a $s$-covering of $sQ$, that is~:
$$ \sup_{x\in sQ} \sharp \left\{j,\ x\in sQ_j  \right\} \lesssim s^d,$$
where $d$ is the homogeneous dimension of the manifold. 
\end{lem}
\dem  We choose $\left( Q(x_j, t^{1/2}/3) \right)_j$ a maximal collection of disjoint balls in $Q$. Then we set $Q_j=Q(x_j, t^{1/2})$, which is a covering of $Q$. 
\\
Fix $x\in sQ$ and denote $J_x:= \left\{j,\ x\in sQ_j  \right\}$. Take  $j_0\in J_x$ (if $J_x\neq \emptyset$ otherwise, there is nothing to prove). By $(D)$, we have
\begin{align*}
\left( \sharp J_x \right) \mu\left(s Q_{j_0}\right) &\lesssim  \left( \sharp J_x \right) s^d \mu\left(\frac{1}{3}Q_{j_0}\right) \\
 & \lesssim  s^d \sum_{j\in J_x} \mu\left(\frac{1}{3}Q_{j}\right) \\
 & \lesssim  s^d  \mu\left(\cup_{j\in J_x} \frac{1}{3}Q_{j}\right) \\
 & \lesssim  s^d  \mu\left(Q(x,2st^{1/2})\right) \\
 & \lesssim  s^d  \mu\left(sQ_{j_0}\right),
\end{align*}
where we  used the fact that the balls $\frac{1}{3}Q_{j}$ are disjoint and have equivalent  measure when the index $j\in J_x$.
\findem

\mb \textbf{Proof of Proposition \ref{prop:pseudo}} Consider a ball $Q$ of radius $r>0$.
We deal with the semigroup and write the oscillation as follows~
$$ f-e^{-r^2\Delta}f = -\int_0^{r^2} \frac{d}{dt} e^{-t\Delta} f dt = \int_0^{r^2} \Delta e^{-t\Delta}f dt.$$
Now we apply arguments used in \cite{ACDH}, Lemma 3.2. Using the completeness of the manifold, we have
\begin{align*}
 \left( \frac{1}{\mu(Q)} \int_Q \left|\int_0^{r^2} \Delta e^{-t\Delta}f dt \right|^q d\mu \right)^{1/q} & \lesssim \int_0^{r^2}\left( \frac{1}{\mu(Q)} \int_{Q} \left| \Delta e^{-t\Delta}f\right|^q d\mu \right)^{1/q}  dt \\
 & \lesssim \int_0^{r^2}\left( \frac{1}{\mu(Q)} \sum_j \int_{Q_j} \left| \Delta e^{-t\Delta}(f-f_{Q_j}) \right|^q d\mu \right)^{1/q}  dt ,
\end{align*}
where $(Q_j)_j$ is a bounded covering of $Q$ with balls of radius $t^{1/2}$ as in Lemma  \ref{cov}.\\
Fix $t\in(0,r^2)$ and denote by $ C_k(Q_j):=2^{k+1}Q_j\setminus 2^kQ_j$ for $k\geq 1$ and $C_0(Q_j)=2Q_j$ . Then, arguing as in Lemma 3.2 in \cite{ACDH}
\begin{align*}
  \sum_j &\int_{Q_j} \left| \Delta e^{-t\Delta}(f-f_{Q_j}) \right|^q d\mu\\ 
  &\lesssim \sum_j \int_{Q_j}  t^{-q}  \left| \int_{M} \frac{ e^{-cd^2(x,y)/t}}{ \mu(Q(y, \sqrt{t}))}(f(y)-f_{Q_j} ) d\mu(y)\right|^q d  \mu(x)
\\&  \lesssim \sum_{j,k; k\geq 0} \int_{Q_j}t^{-q} (\mu(2^{k+1}Q_j))^{q-1} \int_{C_k(Q_j)} \frac{ e^{-cqd^2(x,y)/t}}{ \mu(Q(y, \sqrt{t}))^q}|f(y)-f_{Q_j}|^q d\mu(y) d  \mu(x)
\\
& \lesssim \sum_{j,k; k\geq 1} t^{-q}   (\mu(2^{k+1}Q_j))^{q-1} \int_{C_k(Q_j)}   \left( \int_{\{x;\,d(x,y)\geq 2^{k-1}\sqrt{t}\}} e^{-cqd^2(x,y)/t} d \mu(x)   \right) \frac{|f(y)-f_{Q_j} |^q}{\mu(Q(y, \sqrt{t}))^q} d\mu(y)    
\\  &+\sum_j t^{-q}   \frac{1}{ \mu(Q_j)^q} (\mu(2Q_j))^{q-1} \int_{2Q_j}    \left( \int_{Q_j}d  \mu(x)  \right) \left|f(y)-f_{Q_j} \right|^q d\mu(y)   
\\  & \lesssim \sum_j t^{-q}   \sum_{k \geq 1}  e^{-cq4^k} 2^{kdq}  \int_{C_k(Q_j)}   \left|f(y)-f_{Q_j} \right|^q d\mu(y)    
\\  &+\sum_j t^{-q}   \int_{2Q_j}   \left|f(y)-f_{Q_j} \right|^q d\mu(y)   
\\
 & \lesssim \sum_j t^{-q}   \sum_{k \geq 1}  e^{-cq4^k} 2^{kdq}  \int_{2^{k+1}Q_j}   \left|f(y)-f_{2^{k+1}Q_j} \right|^q d\mu(y)  +\sum_{l=1}^{k+1}\frac{\mu(2^{k+1}Q_j)}{\mu(2^lQ_j)}|f_{2^lQ_j}-f_{2^{l-1}Q_j}|
\\  &+\sum_j t^{-q}   \int_{2Q_j}   \left|f(y)-f_{Q_j} \right|^q d\mu(y)  
\\
& \lesssim \sum_j t^{-q}  \sum_{k\geq 1} e^{-cq4^k} 2^{Mk} t^{q/2} \sum_{l=1}^{k+1} \int_{2^lQ_j} \left|\nabla f \right|^q d\mu+ \sum_j t^{-q}  t^{q/2} \int_{2Q_j} \left|\nabla f \right|^q d\mu.
\end{align*}
We used  (\ref{utp}), $(P_q)$, that for $y\in 2Q_j$,  $\mu(Q(y,\sqrt{t}))\sim \mu(Q_j)$ and for $y\in C_k(Q_j)$, $k\geq 1$, $\frac{1}{\mu(Q(y,\sqrt{t}))}\leq C\frac{2^{kd}}{ \mu(2^{k+1}Q_j)}$. We also used that
for $s$, $t>0$,  
$$\int_{\{x;\,d(x,y)\geq \sqrt{t} \}} e^{-cd^2(x,y)/s} d \mu(x)\leq Ce^{-ct/s}\mu(Q(y, \sqrt{s}))$$ 
thanks to $(D)$ (see Lemma 2.1 in \cite{CD}).

\mb Using that $(2^{l}Q_j)_j$ is a $2^{l}$-bounded covering of $2^{l}Q$, we deduce that
$$ \sum_j \int_{2^{l}Q_j} \left|\nabla f \right|^q d\mu \lesssim 2^{ld} \int_{2^{l}Q} \left|\nabla f \right|^q d\mu \leq 4^{ld}\mu(Q) \sup_{s\geq 1} \aver{sQ} \left|\nabla f \right|^q d\mu,$$
where $d$ is the homogeneous dimension of the doubling manifold.
Thus, it follows that 
\begin{align*}
 \left( \frac{1}{\mu(Q)} \int_Q \left|\int_0^{r^2} \Delta e^{-t\Delta}(f) dt \right|^q d\mu \right)^{1/q} \lesssim \left[\int_0^{r^2} t^{-1/2} dt\right] \sup_{s\geq 1} \left(\aver{sQ} (|\nabla f|^{q})d\mu\right)^{1/q},
\end{align*}
which ends the proof.
\findem

\mb Before we prove off-diagonal estimates under the ``classical'' Poincar\'e inequality, let us recall the following result~:
\begin{prop}\label{RTN} (\cite{AC}) Let $M$ be a complete Riemannian manifold satisfying $(D)$ and $(P_2)$. Then there exists $p_{0}>2$ such that the Riesz transform ${\mathcal R}:=\nabla (-\Delta)^{-\frac{1}{2}}$ is $L^{p}$ bounded for $1<p< p_{0}$.
\end{prop}
We now let 
\be{p0} p_0:=\sup\left\lbrace p\in (2,\infty); \ \nabla (-\Delta)^{-\frac{1}{2}} \textrm { is } L^p \; \textrm {bounded} \, \right\rbrace \tag{$p_0$} \ee 
and 
\be{s0} s_0:=\sup\left\lbrace s\in (1,\infty]; \ (G_s)  \; \textrm {holds} \, \right\rbrace . \tag{$s_0$}  \ee

\begin{rem} \label{rem}  Note that the doubling property $(D)$ and $(DUE)$ imply for $p\in(1,2]$, the $L^p$ boundedness of $ \nabla \Delta^{-\frac{1}{2}}$ which implies $(G_{p})$ (see Subsection \ref{subsec:rapel}) and that $s_0\geq p_0>2$.
\end{rem}  

\mb For the second off-diagonal condition (\ref{off}), we  obtain~:
\begin{prop} \label{prop:ex} Let  $M$ be a complete manifold.
Assume that $M$ satisfies $(D)$ and  admits a classical Poincar\'e inequality $(P_{q})$ for some $q\in[1,\infty)$ as in Definition \ref{classP}. Consider the following estimate
\begin{equation} \label{MH}
 M_{\A,r}(f) \lesssim \M_{q}(|f|+|\nabla f|).
\end{equation} 
\begin{itemize}
\item[1.] If $q_0<2$, then (\ref{MH}) holds 
 for all $r\in(q,s_0)$.
\item[2.] If $q_0\geq 2$, assume moreover $(DUE)$ and that $s_0>q$. Then (\ref{MH}) holds 
 for all $r\in(q,s_0)$.
\end{itemize}
Consequently, 
 (\ref{off}) holds for all $r\in(q,s_0)$.
\end{prop}

\dem It is sufficient  to prove  the following inequalities
\begin{equation} \label{Ue}
\left(\aver{Q}|  e^{-r^2\Delta}f|^rd\mu\right)^{1/r}\leq C \M_{q}(| f |)(x)
\end{equation}
and
\begin{equation} \label{Uge}
\left(\aver{Q}| \nabla e^{-r^2\Delta}f|^rd\mu\right)^{1/r}\leq C \M_{q}(|\nabla f |)(x)
\end{equation}
for every $x\in M$ and every ball $Q$ containing $x$.
We do not detail the proof as it uses analogous argument as in \cite{ACDH}, subsection 3.1,  Lemma 3.2 and the end of this subsection.  For example,  (\ref{Uge}) is essentially inequality (3.12) in section 3 of \cite{ACDH}  where $q_0=2$. 
We just mention that for  (\ref{Ue}), we use the $L^r$ contractivity of the heat semigroup, $(D)$ and $(DUE)$.
 For (\ref{Uge}), we moreover need the following $L^{r}$-Gaffney estimates for $\nabla e^{-t\Delta}$ with $r\in(q_0,s_0)$.  We say that $(\nabla e^{-t\Delta})_{t>0}$ satisfies the $L^p$ Gaffney estimate if there exists $C,\; \alpha>0$ such that for all $t>0$, $E$, $F$ closed subsets of $M$ and $f$ supported in $E$
\begin{equation}\tag{$Ga_p$} \label{gaf}
\|\sqrt{t}|\nabla e^{-t\Delta}f|\|_{L^p(F)}\leq Ce^{-\alpha d(E,F)^2/t}\|f\|_{L^p(E)}.
\end{equation}
In the case where $q_0\geq 2$, interpolating the already known $(Ga_2)$ with $(G_{s})$  for every $ 2<s<s_0$, we get the $(Ga_p)$  for $2<p<s_0$. 
When $q_0< 2$, since in this case  $(G_s)$ holds for all  $1<s<2$ and $2<s<s_0$, interpolating again $(G_s)$ and  $(Ga_2)$,  we obtain the $(Ga_p)$  for all $1<p<s_0$.
\findem
 
\mb It remains to check (\ref{off2}).

\begin{prop} \label{prop1} 
 Let $M$ be a complete manifold satisfying $(D)$ and admitting a classical Poincar\'e inequality $(P_{q})$ for some $1\leq q<\infty$. Then
 \begin{itemize}
\item[1.] If $q_0 < 2$, for  $r>q$, the collection $\A$ satisfies ``$(L^q-L^r)$ off-diagonal'' estimates (\ref{off2}).
\item[2. ]If $q_0\geq  2$, the same result holds under the additional assumption  $(DUE)$.
\end{itemize}
\end{prop}
 
\dem Take $Q_0$, $Q_1$ two equivalent balls, let us say $Q_0 \subset Q_1\subset 10Q_0$ with radius $r_0$ (resp. $r_1$). We choosed a numerical factor $10$ just for convenience. We have to prove that
 \be{amont1}
\left(\frac{1}{\mu(Q_0)}\int_{10Q_0} \left|e^{-r_0^2\Delta}f - e^{-r_{1}^2\Delta}f \right|^r d\mu \right)^{1/r}
\lesssim r_0 \inf_{10Q_0} \M_{q}(|f|+|\nabla f|). \ee
 This is a consequence of 
\be{eqq1} \left(\frac{1}{\mu(Q_0)}\int_{10Q_0}
\left|e^{-r_0^2\Delta}f - e^{-400r_{0}^2\Delta}f \right|^r d\mu \right)^{1/r} \lesssim r_0 \inf_{10Q_0}
\M_{q}(|f|+|\nabla f|)\ee and \be{eqq2} \left(\frac{1}{\mu(Q_0)}\int_{10Q_0} \left|e^{-400r_{0}^2\Delta}f -
e^{-r_{1}^2\Delta}f \right|^r d\mu \right)^{1/r} \lesssim r_0 \inf_{10Q_0} \M_{q}(|f|+|\nabla f|). \ee We use that
$$ e^{-r_0^2\Delta}f - e^{-400r_{0}^2\Delta}f = e^{-r_{0}^2\Delta} \left[1-e^{-399 r_0^2\Delta}\right](f) $$
and $$ e^{-400r_{0}^2\Delta}f - e^{-r_{1}^2\Delta}f = -e^{-r_{1}^2\Delta} \left[1-e^{-(20
r_0)^2-r_{1}^2)\Delta}\right](f).$$ 
We only deal with (\ref{eqq1}), we do the same for (\ref{eqq2}).
From $(D)$ and (\ref{due}), we know that $(UE)$ holds and so we have very fast decays $(L^1-L^\infty)$ for the
semigroup, which permits to gain integrability from $L^q$ to $L^r$. It follows
\begin{align*} 
\lefteqn{\left(\frac{1}{\mu(Q_0)}\int_{10Q_0}
\left|e^{-r_0^2\Delta}f - e^{-400r_{0}^2\Delta}f \right|^r d\mu \right)^{1/r}} & & \\
 & &  \lesssim \sum_{j\geq 0} e^{-\gamma 4^j} \left(\frac{1}{\mu(Q_0)}\int_{C_j(Q_0)}
\left| f-e^{-399 r_0^2\Delta}f \right|^q d\mu \right)^{1/q}, 
\end{align*}
where we make appear the dyadic coronas $C_j(Q_0)$ (see again \cite{ACDH}, Lemma 3.2 and the end of subsection 3.1). Then we use $(D)$ and  $(P_q)$. For each $j$, we choose a bounded covering $(Q_i^j)_i$ of $2^{j+1}Q_0$ with balls of radius $\sqrt{399}r_0$ and obtain
\begin{align*} 
\frac{1}{\mu(Q_0)}\int_{C_j(Q_0)}\left| f-e^{-399 r_0^2\Delta}f \right|^q d\mu  & \lesssim \frac{1}{\mu(Q_0)}\sum_i \int_{Q_i^j}\left| f-e^{-399 r_0^2\Delta}f \right|^q d\mu \\
 & \lesssim \frac{1}{\mu(Q_0)}\sum_i \int_{Q_i^j}\left| f-e^{-399 r_0^2\Delta}f \right|^q d\mu  \\
 & \lesssim \frac{1}{\mu(Q_0)}\sum_i r_0^q \mu(Q_i^j) \sup_{s\geq 1} \left(\aver{sQ_i^j} |\nabla f|^{q}d\mu \right) \\
  & \lesssim \frac{1}{\mu(Q_0)}\sum_i r_0^q \mu(Q_i^j) \sup_{s\geq 1} 2^{dj} \left(\aver{s2^{j+1} Q_0} |\nabla f|^{q}d\mu \right) \\
  & \lesssim \frac{1}{\mu(Q_0)}\sum_i r_0^q \mu(Q_i^j) 2^{dj} \inf_{Q_0} \M \left(|\nabla f|^{q}\right)  \\
  & \lesssim r_0^q  2^{dj} \frac{\mu(2^{j+1}Q_0)}{\mu(Q_0)}  \left( \inf_{Q_0} \M_{q}(|\nabla f|)\right)^q  \\
  & \lesssim r_0^q  2^{2dj}  \left(\inf_{Q_0} \M_{q}(|\nabla f|)\right)^q.
\end{align*}
We applied $(P_q)$ in the third inequality. In the fourth inequality, we used that $sQ_i^j\subset 2^{j+1}s Q_0$ and  thanks to $(D)$, $\mu(2^{j+1} s Q_0) \lesssim \mu(sQ_i^j) 2^{jd}$. Then we applied the bounded overlap property in the sixth one. \\
Summing in $j$, we show the desired inequality (\ref{eqq1}). Similarly we prove (\ref{eqq2}), which completes the proof of (\ref{amont1}). 
\findem

\mb We get the following corollary:

\begin{cor} \label{corfinal} Assume that $M$ is complete, satisfies $(D)$ and admits a classical Poincar\'e inequality $(P_q)$ for some $q\in[1,\infty)$. In the case where  $q_0\geq 2$, we moreover assume $(DUE)$ and $s_0>q$. Then the  assumptions of Theorem \ref{thm:CZ} and \ref{thm:int} hold. We have pseudo-Poincar\'e inequality $(\widetilde{P}_{q})$ and $\A$ satisfies ``$L^q-L^r$ off-diagonal estimates''  for $ r\in (q,s_0)$. 
\end{cor}
 
\mb {\bf Conclusion : } When $q<2$, the assumptions of Theorem \ref{thm:CZ} (according to this particular choice of $\A$) are weaker than the Poincar\'e inequality and are sufficient to get the Calder\'on-Zygmund decomposition.

\mb We also have the homogeneous version: 
  \begin{cor} \label {ISPH} Assume that $M$ is complete, satisfies $(D)$ and admits a classical Poincar\'e inequality $(P_q)$ for some $1\leq q<\infty$. In the case where  $q_0\geq 2$, we moreover assume $(DUE)$. \\
 Let $\A:=(A_Q)_Q$ with $A_Q:=e^{-r_Q^2\Delta}$. Then the  assumptions of Theorems \ref{thm:CZH} and \ref{ISH} holds. We have  pseudo-Poincar\'e inequality $(\widetilde{P}_{q})$, $\A$ satisfies ``homogeneous $L^q-L^r$ off-diagonal estimates''  for $ r\in (q,s_0)$. 
\end{cor}

\subsection{Application to Reverse Riesz transform inequalities.}

We refer the reader to \cite{AC,ACDH} for the study of the so-called (\ref{rrp}) inequalities~:
\begin{equation} \label{rrp}
 \| \Delta^{1/2} f\|_{L^p} \lesssim \| |\nabla f| \|_{L^p}.  \tag{$RR_p$}
\end{equation}

\mb We know that $(RR_2)$ is always satisfied and that $(D)$ and $(DUE)$ implies (\ref{rrp}) for all $p\in(2,\infty)$. For the exponents lower than 2, P. Auscher and T. Coulhon obtained the following result (\cite{AC})~:

\begin{thm} Let $M$ be a complete non-compact doubling Riemannian manifold. Moreover assume that the classical Poincar\'e inequality $(P_q)$ holds for some $q\in (1,2)$.  Then for all $p\in(q,2)$, (\ref{rrp}) is satisfied.
\end{thm}

\mb This result is based on a Calder\'on-Zygmund decomposition for Sobolev functions. Using our new assumptions, we also obtain the following improvement~:

\begin{thm} Assume that $M$ is complete, satisfies $(D)$ and admits a pseudo-Poincar\'e inequality $(\widetilde{P}_{q})$ for some $q\in(1,2)$. If in addition, the collection $\A$ satisfies $L^q-L^2$ ``off-diagonal estimates'', 
 then  (\ref{rrp}) holds for all $p\in(q,2)$.
\end{thm}

\begin{rem} Corollary \ref{corfinal} shows that these new assumptions are weaker than the Poincar\'e inequality $(P_{q})$.
\end{rem}

\mb We do not prove this result and refer the reader to \cite{AC}. The proof is exactly the same as it relies on the Calder\'on-Zygmund decomposition.

\begin{rem} We refer the reader to other works of the authors \cite{max,BB}. In \cite{max}, the assumption (\ref{rrp}) plays an important role in order to prove some maximal inequalities in dual Sobolev spaces $W^{-1,p}$, which do not require Poincar\'e inequalities. So it might be important to know how to  prove (\ref{rrp}) without Poincar\'e inequality.
\end{rem}

\subsection{Application to Gagliardo-Nirenberg inequalities.}

We devote this subsection to the study of Gagliardo-Nirenberg inequalities. We refer the reader to \cite{Nadine2} for a recent work on this subject.

\begin{df} We introduce the Besov space. For $\alpha<0$, we set $B_{\infty,\infty}^{\alpha}$ the set of all measurable functions $f$ such that
$$
\|f\|_{B_{\infty,\infty}^{\alpha}}:=\sup\limits_{t>0}t^{-\frac{\alpha}{2}}\|e^{-t\Delta}f\|_{L^\infty}<\infty.
$$
We have the following equivalence (Lemma 2.1 in \cite{Nadine2})~:
$$ \|f\|_{B_{\infty,\infty}^{\alpha}}\sim\sup_{t>0}t^{-\frac{\alpha}{2}}\|e^{-t\Delta}(f-e^{-t\Delta}f)\|_{L^\infty}. $$
\end{df}

\mb Then, the so-called Gagliardo-Nirenberg inequalities are~:
\begin{equation}\label{led}
\|f\|_{l}\lesssim  \|\,|\nabla f|\,\|_{p}^{\theta} \|f\|_{B_{\infty,\infty}^{\frac{\theta}{\theta-1}}}^{1-\theta}
\end{equation}
 where $\theta= \frac{p}{l}$ for some $p,l\in[1,\infty)$.

\mb We first recall one of the main results of \cite{Nadine2}:

\begin{thm}\label{G} Let $M$ be a complete non-compact Riemannian manifold satisfying $(D)$ and $(P_{q})$  for some $1\leq q<\infty$. Moreover, assume that $M$ satisfies the global pseudo-Poincar\'{e} inequalities $(P'_{q})$ and $(P'_{\infty})$. Then (\ref{led}) holds for all $q\leq p<l<\infty$.
\end{thm}

\mb Here, the \textbf{global} pseudo-Poincar\'e inequality $(P'_{q})$ for some $q\in[1,\infty]$ corresponds to
\begin{equation*}\tag{$P'_{q}$}
\|f-e^{-t\Delta}f\|_{L^q}\leq Ct^{\frac{1}{2}}\|\,|\nabla f|\,\|_{L^q}.
\end{equation*} 

\mb This result requires global pseudo-Poincar\'e inequalities and some Poincar\'e inequalities with respect to balls. These two kinds of inequalities are quite different as they deal with oscillations with respect to the semigroup (for the pseudo-Poincar\'e inequalities) and to the mean value operators (for the Poincar\'e inequalities). We saw in the previous subsection, that Poincar\'e inequality implies pseudo-Poincar\'e inequality. That is why, we are looking for assumptions requiring only the Poincar\'e inequality, getting around the assumed global pseudo-Poincar\'e inequalities.

\mb We begin first showing that pseudo-Poincar\'e inequalities related to balls yield global pseudo-Poincar\'e inequalities.
 
 \begin{prop} \label{locglob} Let $M$ be a complete Riemannian manifold satisfying $(D)$ and admitting a pseudo-Poincar\'e inequality $(\widetilde{P}_q)$ for some $1\leq q<\infty$. Then the global pseudo-Poincar\'e inequality $(P'_{q})$ holds.
\end{prop}

\dem Let $t>0$. Pick a countable set
$\left\lbrace x_{j}\right\rbrace _{j\in J} \subset  M,$ such that $ M=
\underset{j\in J}{\bigcup}Q(x_{j},\sqrt{t}):=\underset{j\in J}{\bigcup}Q_{j}$ 
and for all $x\in M$,  $x$ does not belong to more than $N_{1}$ balls $Q_{j}$.
Then
 \begin{align*}
 \|f-e^{-t\Delta}f\|_q^q  &\leq \sum_{j}\int_{Q_j}|f-e^{-t\Delta}f|^q d\mu \\
 & \lesssim \sum_{j}  t^{\frac{q}{2}} \int_{Q_j}|\nabla f|^qd\mu \\
 & \lesssim N_1 t^{\frac{q}{2}} \int_M |\nabla f|^qd\mu.
\end{align*}
\findem

\begin{rem} It is easy to see that the global pseudo-Poincar\'e inequality $(P'_{\infty})$ is satisfied under  $(D)$ and $(DUE)$ (see for instance \cite{Nadine2}, p.499).
\end{rem}

\mb Using Propositions  \ref{locglob}, \ref{prop:pseudo} and Theorem \ref{G}, we get the following improvement version of Theorem 1.2 in \cite{Nadine2}~:

\begin{thm} Let $M$ be a complete Riemannian manifold satisfying $(D)$ and admitting a Poincar\'e inequality $(P_q)$ for some $1\leq q <\infty$. If $ q_0\geq 2$, we moreover assume $(DUE)$. Then (\ref{led}) holds for all $q\leq p<l<\infty$.
 \end{thm}

\mb Using our new assumptions, we get also the following Gagliardo-Nirenberg theorem:

\begin{thm}\label{gag2} Assume that $M$ satisfies the hypotheses of Theorem \ref{ISH} with $A_Q=e^{-r_Q^2\Delta}$ and that $r=\infty$. Moreover, we assume $(DUE)$. Then (\ref{led}) holds for all $q\leq p<l<\infty$.
\end{thm}
\dem
The proof is analogous to  that of Theorems 1.1 and 1.2 in \cite{Nadine2}. We use our homogeneous interpolation result of Theorem \ref{ISH}. Also we need our non-homogeneous interpolation result of Theorem \ref{thm:int}. It holds thanks to (\ref{Ue}) which is true under $(D)$ and $(DUE)$. Moreover, $(P'_q)$ is satisfied and $(P'_{\infty})$ holds thanks to $(D)$ and $(DUE)$. 
\findem

\mb As a Corollary, we obtain
\begin{thm}\label{lp} Consider a complete Riemannian manifold $M$ satisfying $(D)$, $(P_{q})$ for some $1\leq q <\infty$ and assume that there exists $C>0$ such that for every $x,\,y\in M$ and $t>0$
\begin{equation*}\tag{$G$}
|\nabla_{x}p_{t}(x,y)|\leq \frac{C}{\sqrt{t}\mu(B(y,\sqrt{t}))}.
\end{equation*}
($(G)$ is equivalent to the assumption $(G_{\infty})$.)  In the case where $q_0>2$, we moreover assume $(DUE)$. 
Then inequality (\ref{led}) holds for all $q \leq p<l<\infty$.
\end{thm}
\dem In the case where $q\leq 2$, this result is already in \cite{Nadine2}.  For $q_0\geq 2$,  we are under the hypotheses of Theorem \ref{gag2} thanks to subsection 4.1 and since $(G)$ implies that $r=\infty$.
\findem


\begin{thebibliography}{100}

\bibitem{ambrosio1}
L.~Ambrosio, M.~Miranda~Jr, and D.~Pallara.
\newblock Special functions of bounded variation in doubling metric measure
  spaces.
\newblock {\em Calculus of variations~: topics from the mathematical heritage of
  E. De Giorgi, Quad. Mat., Dept. Math, Seconda Univ. Napoli, Caserta},
  {\textbf 14}, pages 1--45, 2004.

\bibitem{A1}
P.~Auscher.
\newblock On ${L}^p$ estimates for square roots of second order elliptic
  operators on ${\R}^{n}$.
\newblock {\em Publ. Mat. \textbf{48}}, pages 159--186, 2004.

\bibitem{A}
P.~Auscher.
\newblock On necessary and sufficient conditions for ${L}^p$ estimates of
  {R}iesz transforms associated to elliptic operators on ${\R}^n$ and related
  estimates.
\newblock {\em Memoirs of Amer. Math. Soc. \textbf{186} no.871}, 2007.

\bibitem{A2}
P.~Auscher.
\newblock On the Calder\'on-Zygmund lemma for Sobolev functions: http://arxiv.org/abs/0810.5029.

\bibitem{AB}
P.~Auscher and B.~Ben Ali.
\newblock Maximal inequalities and {R}iesz transform estimates on ${L^{p}}$ spaces for {S}chr\"{o}dinger operators  with nonnegative potentials.
\newblock {\em  Ann. Inst. Fourier, \textbf{57}, no 6}, pages 1975--2013,  2007.

\bibitem{AC}
P.~Auscher and T.~Coulhon.
\newblock Riesz transform on manifolds and Poincar\'e inequalities.
\newblock {\em Ann. Sc. Nor. Sup. Pisa (\textbf{5}), IV, 3}, pages 531--555,
  2005.

\bibitem{ACDH}
P. Auscher, T. Coulhon, X. T. Duong and S. Hofmann.
  \newblock Riesz transform on manifolds and heat kernel regularity,
     \newblock {\em Ann. Sci. Ecole Norm. Sup. 37}, pages 911-957, 2004.



\bibitem{AM}
P.~Auscher and J.M.~ Martell.
\newblock Weighted norm inequalities, off-diagonal estimates and elliptic
  operators. {Part I} : {G}eneral operator theory and weights.
\newblock {\em Adv. in Math. \textbf{212}}, pages 225--276, 2007.




\bibitem{Nadine5}
N.~Badr.
\newblock Ph.{D} {T}hesis.
\newblock {\em  Universit\'e Paris-sud 11}, 2007.

\bibitem{Nadine}
N.  Badr.
\newblock Real interpolation of Sobolev spaces,
\newblock {\em  to appear in Math. Scand.}

\bibitem{Nadine1}
N.~Badr.
\newblock Real interpolation of Sobolev spaces associated to a weight,
\newblock {\em to appear in Pot. Anal.}

\bibitem{Nadine2}
N.~Badr.
\newblock Gagliardo-Nirenberg inequalities on manifolds.
\newblock {\em J.M.A.A, {\textbf 349}}, pages 493--502, 2009.

\bibitem{Nadine3}
N.~Badr and E.~ Russ.
\newblock Interpolation of {S}obolev Spaces, {L}ittlewood-{P}aley inequalities and {R}iesz transforms on graphs,
\newblock {\em Pub. Mat., {\textbf 53 (2)}}, pages 273--328, 2009.

\bibitem{Nadine4}
N.~Badr and B.~Ben Ali.
\newblock ${L}^p$ boundedness of {R}iesz tranform related to {S}chr\"odinger operators on a manifold,
\newblock {\em  to appear in Ann. Sc. Nor. Sup. Pisa.}

\bibitem{BB}
N.~Badr and F.~ Bernicot.
\newblock Abstract Hardy-Sobolev spaces and interpolation,
\newblock {\em  submitted}, 2009 available at http://arxiv.org/abs/0901.0518.


 \bibitem{Besma}
B.~Ben Ali.
\newblock Ph.{D} {T}hesis.
\newblock {\em  Universit\'e Paris-sud 11}, 2008.

\bibitem{bennett}
C.~Bennett and R.~Sharpley.
\newblock {\em Interpolation of operators}.
\newblock Academic Press, 1988.

\bibitem{bergh}
J.~Bergh and J.~L\"{o}fstr\"{o}m.
\newblock {\em Interpolations spaces, {A}n introduction}.
\newblock Springer (Berlin), 1976.

\bibitem{BJ}
F.~Bernicot and J.~Zhao.
\newblock Abstract {H}ardy spaces, {\em J. Funct. Anal.
  \textbf{255}}, no.~7, pages 1761--1796, 2008.

\bibitem{B2}
F.~Bernicot.
\newblock Use of abstract {H}ardy spaces, real interpolation and
  applications to bilinear operators, {\em Mathematische Zeitschrift} (2009)
  available at http://fr.arxiv.org/abs/0809.4110.

\bibitem{max}
F.~Bernicot.
\newblock Maximal inequalities for dual Sobolev spaces $W^{-1,p}$ and applications to interpolation,
\newblock {\em Mathematical Research Letters}, 2009.

\bibitem{coifman2}
R.~Coifman and G.~Weiss.
\newblock {\em Analyse harmonique sur certains espaces homog\`{e}nes}.
\newblock Lecture notes in Math., Springer, 1971.

\bibitem{CD1}
T.~Coulhon and X.T. Duong.
\newblock Riesz transforms for $1\leq p\leq 2$.
\newblock {\em Trans. Amer. Math. Soc. \textbf{351} no.2}, pages 1151--1169,
  1999.

\bibitem{CD}
T.~Coulhon and X.T. Duong.
\newblock Maximal regularity and kernel bounds: observations on a theorem by
  {H}ieber and {P}r\"uss.
\newblock {\em Adv. Differential Equations \textbf{5} no.1-3}, pages 343--368,
  2000.

\bibitem{Davies2}
E.~B. Davies.
\newblock Non-{G}aussian aspects of heat kernel behavior.
\newblock {\em J. London Math. Soc.}, 1997.

\bibitem{DY}
X.T. Duong, L.~Yan.
\newblock Duality of {H}ardy and {BMO} spaces associated with operators with
  heat kernel bounds.
\newblock {\em J. Amer. Math. Soc. \textbf{18}}, no.4, pages 943--973, (2005).

\bibitem{goldshtein}
  V. Gol'dshtein and  M. Troyanov.
  \newblock   {A}xiomatic {T}heory of {S}obolev {S}paces.
  \newblock {\em Expo. Mathe., \textbf{19}}, pages 283-336, 2001.
 
\bibitem{gri}
A.~Grigor'yan.
\newblock Gaussian upper bounds for the heat kernel on arbitrary manifolds.
\newblock {\em J. Diff. Geom.} \textbf{45}, pages 33--52, 1997.

\bibitem{hajlasz4}
P.~Hajlasz and P.~Koskela.
\newblock {S}obolev met {P}oincar\'{e}.
\newblock {\em Mem. Amer. Math. Soc. \textbf{145}} no. 688, pages 1--101, 2000.

\bibitem{HM}
S.~Hofmann and S.~Mayboroda.
\newblock Hardy and {BMO} spaces associated to divergence form elliptic
  operators.
\newblock {\em Math. Ann., to appear}.

\bibitem{KZ}
S.~Keith and X.~Zhong.
\newblock The Poincar\'e inequality is an open ended condition.
\newblock {\em Ann. of Maths. \textbf{167} no. 2}, pages 575--599, 2008.

\bibitem{LY}
P. Li and S.T. Yau, 
\newblock On the parabolic kernel of the {S}chr{\"o}dinger operator. 
\newblock {\em Acta Math. \textbf{156}}, pages 153--201, 1986. 

\bibitem{Martell}
J.M.~Martell.
\newblock Sharp maximal functions associated with approximations of the identity in spaces of homogeneous type and applications.
\newblock {\em Studia Math. \textbf{161}}, pages 113--145 (2004).

\bibitem{saloff}
L. Saloff-Coste. 
\newblock A note on Poincar\'e, Sobolev and Harnack inequalities. 
\newblock{\em Duke J. Math. \textbf{65}}, pages 27--38, 1992. 

\bibitem{Se}
E.M. Stein,
 Singular integrals and differentiability properties of
  functions,
 Princeton Univ. Press, 1970.

\end{thebibliography}
\end{document}